\documentclass[reqno,12pt]{amsart}

\usepackage{amsmath}
\usepackage{amsfonts}
\usepackage{amssymb}
\usepackage{eufrak}
\usepackage{amscd}
\usepackage{amsthm}
\usepackage{epsfig}
\usepackage{amstext}
\usepackage[all]{xy}

\theoremstyle{plain}

 \newtheorem{theorem}{Theorem}

 \newtheorem{cor}[theorem]{Corollary}
 \newtheorem{lemma} [theorem] {Lemma}
 \newtheorem{remark}[theorem]{Remark}

\begin{document}

 \title{Semiprojectivity of  universal C*-algebras generated by algebraic elements}
\author{Tatiana Shulman}

\address{Department of Mathematical Sciences, University of Copenhagen, Universitetsparken 5, DK-2100 Copenhagen, Denmark}

\email{shulman@math.ku.dk}



\date{\today}

\maketitle

\begin{abstract} Let $p$ be a polynomial in one variable whose  roots either all have multiplicity more than 1 or all have multiplicity exactly 1. It is shown that the universal $C^*$-algebra of a relation $p(x)=0$, $\|x\| \le 1$ is semiprojective. In the case of all roots multiple it is shown that the universal $C^*$-algebra is also residually finite-dimensional. Applications to polynomially compact operators are given.
\end{abstract}

\section*{Introduction}
The "general lifting problem" in C*-algebras is a class of problems which are formulated in the following way:
for a quotient $C^*$-algebra $A/I$ and its element $b$  with some specific properties,  to find (or to prove the
existence of) a preimage of $b$ in $A$ with the same properties.
Lifting problems are closely related to the theory of stability of relations which,  roughly speaking,
investigates when "an almost representation of a relation" is "close" to a representation (for strict definitions  see \cite{Loring}).

 By a relation we mean a system of equalities of
the form $f(x_1,...,x_n,x_1^*,...,x_n^*) = 0$, where $f$ is a non-commutative polynomial (an element of the
corresponding free algebra), which can include also inequalities of the form $\|x_i\|\le c$ or $\|x_i\|<c$ (they also can be written as polynomial equalities if one extends the list of generators).

Recall that a $C^*$-algebra $D$ is  {\it projective} if for any $C^*$-algebra $A$, its ideal $I$ and every
$\ast$-homomorphism $\phi:D \to A/I$, there exists  $\tilde \phi$ such that the diagram
$$\xymatrix{ & A \ar[d] \\ D\ar[ru]^{\tilde \phi}\ar[r]^{\phi} & A/I}$$ commutes.

A $C^*$-algebra $D$ is  {\it semiprojective}  if for any $C^*$-algebra $A$, any increasing chain of ideals $I_n$
in $A$  and for every $\ast$-homomorphism $\phi:D \to A/I$, where   $I=\overline{\bigcup I_n}$, there exist
 $n$ and  $\ast$-homomorphism $\tilde \phi : D\to A/I_n$ such that the diagram
$$\xymatrix{ & A/I_n \ar[d] \\ D\ar[ru]^{\tilde \phi}\ar[r]^{\phi} & A/I}$$ commutes.

A $C^*$-algebra $D$ is  {\it weakly semiprojective}  if for any
sequence  of $C^*$-algebras $A_i$ and any $\ast$-homomorphism
$\phi:D \to \prod A_i /\bigoplus A_i$, there exists a
$\ast$-homomorphism  $\tilde \phi: D \to \prod A_i$ such that the
diagram
$$\xymatrix{ & \prod A_i \ar[d] \\ D\ar[ru]^{\tilde \phi}\ar[r]^{\phi} &\prod A_i /\bigoplus A_i}$$ commutes.

These notions provide an algebraic setting of lifting and stability questions: a relation is liftable (stable,
weakly stable) if and only if its universal $C^*$-algebra is projective (semiprojective, weakly semiprojective),
\cite{Blackadar}, \cite{Loring}. Of course it is assumed that the universal $C^*$-algebra of a relation exists.
 Note that  in general a relation needn't have a universal $C^*$-algebra, but if it includes restrictions
$\|x_i\|\le c_i$, for all generators,  then the universal $C^*$-algebra exists (of course if the relation has
at least one representation in a Hilbert space).

Thus  given a liftable or (weakly) stable relation there arises a natural question if being combined with
restrictions on norms of generators it is still liftable or (weakly) stable: the positive answer would produce
new examples of projective or (weakly) semiprojective $C^*$-algebras.

It was proved in  \cite{Olsen} that each nilpotent is liftable and
 in \cite{Loring}, \cite{LorPed} it
was asked a question if a relation
\begin{equation*}\label{rel} x^n=0,  \|x\| \le 1
\end{equation*} is
liftable or, in other words, if  it defines a  projective $C^*$-algebra. In \cite{nilpotents}  this question was
answered positively and the main tool was using the technique of M-ideals in Banach spaces.

In \cite{Don} it was proved that for each polynomial $p$ in one
variable,  a relation $p(x)=0$ is weakly stable (though the proof
given in \cite{Don} easily implied also stability) but again it
remained unknown if the  universal $C^*$-algebra of a relation
$$p(x)=0, \|x\| \le 1$$ is semiprojective. In the present paper using  approach of \cite{nilpotents} we
prove this in two "opposite" cases:   when all the roots of the polynomial have multiplicity more than 1 and when all the roots have multiplicity 1.

We prove also that for polynomials with all roots multiple, the universal $C^*$-algebras of such relations are residually finite-dimensional (have separating families of
finite-dimensional representations). Note that this in no way means that these algebras are small. For
polynomials of degree more than two they all are non-exact.

Theorems \ref{Calkin} and \ref{Calkin2} concern well known questions of C. Olsen. Let $T\in B(H)$, $p(x)$ be a polynomial in one
variable. Olsen's questions are:

1) Does there exist a compact operator $K$ such that $\|p(T+K)\|=\|p(T)\|_e$?

\noindent(by $\|T\|_e$  the essential norm of $T$ is denoted).

2) Can this $K$ be chosen common for all polynomials?

The questions are still open but there are partial results (\cite{Akemann}, \cite{OlsenAlg}, \cite{OlsenPlast}, \cite{OlsenPerturb},
 \cite{SmithAndK}, \cite{SmithWard}). In particular Olsen proved (\cite{OlsenAlg}) that under the assumption that
$p(T)$ is compact, such $K$ always exists. In this paper we prove that under the same assumption $K$ can be chosen common for polynomials $p(x)$ and $x$ in the case when either all roots of the polynomial  have multiplicity more than 1 or when they all have multiplicity 1. In other words we prove (Theorem \ref{Calkin}) that for such polynomials if $p(T)$ is compact then there exists a compact operator $K$ such that $$p(T+K)=0, \;\;\; \|T+K\| = \|T\|_e.$$

\section{Preliminary results and constructions}

Below we use the following notation.

By $M(I)$ we denote the multiplier $C^*$-algebra of a
$C^*$-algebra $I$.

For a closed ideal $I$ of a $C^*$-algebra $A$, we denote by $a\to \dot a$ the standard epimorphism from $A$ to $A/I$. We say
that $a$ is a lift of $b$ if $\dot{a} = b$.

Given ideals $I_1\subseteq I_2\subseteq \ldots \subseteq
I\subseteq A$ with $I=\overline{\bigcup I_n}$ we use the notation
$$q_n: A\to A/I_n$$ and $$q_{n, \infty}:A/I_n \to A/I$$ for the
induced epimorphisms.

\medskip

For an arbitrary algebra $A$, an operator $T:A\to A$ is
called an elementary operator if it is of the form $T =
\sum_{i=1}^N L_{a_i}R_{b_i}$, where $a_i, b_i\in A^+$. By $L_{a_i}, R_{b_i}$
we denote operators of left and right multiplication by $a_i$ and $b_i$ respectively.

A subspace $Y$ of a Banach space $X$ is called proximinal  if for
any $x\in X$ there is $y\in Y$ such that $\|x-y\| = d(x, Y)$, or,
equivalently, if for any $z\in X/Y$ there is a lift $x$ of $z$
such that $\|z\|=\|x\|$.

Below for any $C^*$-algebra $A$ we denote the set of its self-adjoint elements by $A_{sa}$.

\begin{theorem}\label{prox}(\cite{nilpotents})
Let $A$ be a $C^*$-algebra, $T$ an elementary operator on $A$.
Then, for each ideal $I$ of $A$,  $\overline{TI}$ is
proximinal in $\overline{TA}$ and  $\overline{TI_{sa}}$ is
proximinal in $\overline{TA_{sa}}$.
\end{theorem}

\begin{remark} In \cite{nilpotents} it was not considered the case of $\overline{TI_{sa}}$ but the proof is the same, if one considers
 $C^*$-algebras as Banach spaces over $\mathbb R$.
 \end{remark}

Theorem \ref{prox} implies a corollary which was contained in nondirect form in
the proof of [Theorem 5,  \cite{nilpotents}]:

\begin{cor}\label{new} Let $A$ be a $C^*$-algebra, $I$  its  ideal, $T$ an elementary operator on $A$, $b\in A/I$.
 If, for some lift $a_0$ of $b$, an element $Ta_0$ is also a lift of $b$, then
 the set $\overline{\{ Ta| \; \dot a = b\}}$ contains an element
 of norm equal to $\|b\|$.
\end{cor}
\begin{proof}
Let $E= \{ Ta| \; \dot a = b\}$.  Then
$$E = \{Ta_0 +Ti\;|\;\;i\in I\},\;\; \overline E = \{Ta_0
+x\;|\;\;x\in \overline{TI}\}.$$

All elements of $\overline E$ are lifts of $b$. Indeed let $\dot a
= b$, then $a=a_0+i$, for some $i\in I$. Since $T$ is elementary,
it preserves ideals, whence $(Ti)^{\cdot}=0$ and $(Ta)^{\cdot} =
(Ta_0)^{\cdot} = b$. This easily implies that $\overline E$ also consists of lifts of
$b$.

Let $b_i, c_i, i=1, \ldots, N$, be elements of $A^+$ such that $A=\sum_{i=1}^N L_{b_i}R_{c_i}$ and  let $B$
be the  $C^*$-algebra generated by $b_i,  i=1,\ldots, N$, and $I$. Let us choose  a quasicentral approximate
unit $u_{\lambda}$ for the ideal $I$ of $B$. Then,  for each $a\in E$, we have
$T\left((1-u_{\lambda})a\right)\in E$ and $$\lim \|T\left((1-u_{\lambda})a\right)\| = \lim \|(1-u_{\lambda})Ta\|
= \|(Ta)^{\cdot}\|=\|b\|$$ (we used here a well known fact \cite{Murphy}, Section 3.1.3, that if
$e_{\lambda}$ is any approximate unit then
 $$\lim_{\lambda} \|x(1-e_{\lambda})\| = \|\dot
x\|,$$ for any $x\in A$). Hence $\inf_{y\in \overline{TI}} \|Ta_0
+y\| \le \|b\|$. Since the norm of any lift of $b$ is not less
than $\|b\|$ we get $\inf_{y\in \overline{TI}} \|Ta_0 +y\| =
\|b\|$.

Thus $dist(Ta_0, \overline{TI}) = \inf_{y\in \overline{TI}} \|Ta_0
+y\| = \|b\|$ and, by Theorem \ref{prox}, there is $y_0\in
\overline{TI}$ such that $\|Ta_0+y_0\| = \|b\|$. Since
$Ta_0+y_0\in \overline E$, we are done.
\end{proof}

 \begin{lemma}\label{aAa} $\overline{\dot a (A/I) \dot a} =
 \overline{aAa}/\overline{aIa}$.
 \end{lemma}
 \begin{proof} Define a map $f: aAa \to \dot a (A/I) \dot a$ by
 $f(aba)=\dot a \dot b \dot a$ and extend it by continuity to
 $\tilde f: \overline{aAa}\to \overline{\dot a (A/I) \dot a}$.
 Since $\tilde f|_{\overline{aIa}}=0$ there is a well-defined map
 $F: \overline{aAa}/\overline{aIa} \to \overline{\dot a (A/I) \dot
 a}$. This is surjection because its image is closed (as the image of any $\ast$-homomorphism between $C^*$-algebras) and contains $\dot a (A/I) \dot
 a$.

 Suppose that $x\in \overline{aAa}/\overline{aIa}$ and
 $F(x)=0$. Let $X\in \overline{aAa}$ be any preimage of $x$. Then
 there are elements $b_n\in A$ such that $X=\lim ab_na$. Then
 $\lim \dot a\dot {b_n}\dot a=0$ and there are $i_n\in I$ such
 that $ab_na-i_n\to 0$. Hence $X\in I$. Let $\{u_{\lambda}\}$ be a
 quasicentral approximate unit for $I\subseteq A$. Then  for any $\epsilon$, there
 is
  $n$ such that $\|X-ab_na\|\le \epsilon$ and there is $\lambda $
  such that $\|X-Xu_{\lambda}\|\le \epsilon$,
  $\|ab_nau_{\lambda}-ab_nu_{\lambda}a\|\le \epsilon$. Then
  $\|X-ab_nu_{\lambda}a\|\le 3\epsilon$, $dist(X, aIa)\le
  3\epsilon$ whence $X\in \overline{aIa}$.
 \end{proof}

 A $C^*$-algebra $A$ is called SAW* (\cite{Pedersen}) if for any $a, b\in A_{sa}$,
 $ab=0$,
 there is $0\le z\le 1$, $z\in A^+$, such that $az=a$, $zb=0$.

 \begin{lemma}\label{SAW} If $A$ is SAW*, $p\in A$ is a projection, then
 $pAp$ is SAW*.
 \end{lemma}
 \begin{proof} Suppose $a, b\in (pAp)_{sa}$, $ab=0$. Since $A$ is
 SAW* there is $0\le z\le 1$, $z\in A^+$, such that $az=a$,
 $zb=0$. Clearly $0\le pzp\le 1$, $pzp\in (pAp)^+$,  $apzp=a$,
 $pzpb=0$.
 \end{proof}

Below we write $d<<e$ for selfadjoint elements satisfying $de = ed
=d$.

 \begin{theorem}\label{Olsen, Pedersen} (\cite{Olsen}) Let $A$ be a $C^*$-algebra,
 $I$ --- its ideal, such that the quotient $A/I$ is SAW*. Let $b\in
 A/I$ be a nilpotent of order n,  $n\ge 2$. Then there are elements \begin{equation}\label{2}0\le
 d_{n-1}<<e_{n-1}<<d_{n-2}<<e_{n-2}<<\ldots <<e_2<<d_1<<e_1\le
 e_0=1\end{equation} and their lifts \begin{equation}\label{3}0\le
 D_{n-1}<<E_{n-1}<<D_{n-2}<<E_{n-2}<<\ldots <<E_2<<D_1<<E_1\le
 E_0=1\end{equation} such that

 (i) $b = \sum_{j=1}^{n-1}(e_{j-1}-e_j)bd_j$,

 (ii) for any lift $B$ of $b$, $\sum_{j=1}^{n-1}(E_{j-1}-E_j)BD_j$
 is a nilpotent lift of $b$ of order $n$.
 \end{theorem}

Now let an element $b$ of a $C^*$-algebra  be  algebraic  with the
minimal polynomial $p(t) =
 (t-\lambda_1)^{k_1}\ldots (t-\lambda_n)^{k_n}$. In \cite{Don} a family of
 $n$ pairwisely orthogonal projections was associated with $b$ in the following way.
  Let $e_i$ be the spectral idempotents for
 $b$ and $s= \sum e_i^*e_i$. One can prove that $e_i(b-\lambda_i)^{k_i}e_i=0$ for each $i$,  $ s$ is
 invertible and \begin{equation} \label{projections}p_i = s^{1/2}e_is^{-1/2}\end{equation}
  are
 projections with sum 1. Let $c=s^{1/2}bs^{-1/2}.$ Then
 \begin{equation} \label{1}
 p_i(c-\lambda_i)^{k_i}p_i=0,\end{equation}

 \begin{equation} \label{4}
\sum_{i=1}^n p_icp_i=c.\end{equation}

\begin{theorem}(\cite{Don})  The following are equivalent:
\item (i) $b\in A/I$ has a lift
$\tilde b\in A$ with $p(\tilde b)=0$;

\item(ii) the family $\{p_i\}$ can be lifted
to a family of
 projections in $A$ with sum 1.
 \end{theorem}
 Since in what follows we use ideas from \cite{Don}, we now sketch arguments from \cite{Don} of how (ii) implies (i).
 Let a family of projections $\{Q_i\}$ with $\sum
Q_i =1$ be a lift of $\{p_i\}$. Since $p_i(c-\lambda_i)p_i$ is a nilpotent of order $k_i$
in a $C^*$-algebra $p_iA/Ip_i = Q_iAQ_i/Q_iIQ_i$, by \cite{Olsen}
it has a nilpotent lift $t_i\in Q_iAQ_i$ of order $k_i$. Now one
can check that $\sum t_i+\lambda_iQ_i$ is a lift of $c$ and
$p(\sum t_i+\lambda_iQ_i)=0$. As a lift of $b$ we can take
$S^{-1/2}(\sum t_i+\lambda_iQ_i)S^{1/2}$, where $S$ is any
invertible positive lift of $s$ (to find such $S$ we can take any
selfadjoint lift of $s$ and then use functional calculus).

\section{Main results}

\begin{lemma}\label{main} Let $p(t) ~=
(t-\lambda_1)^{k_1}\ldots (t-\lambda_n)^{k_n}$ be a polynomial such that $k_i\ge 2$, $i=1, \ldots, n$, $A$ be a $C^*$-algebra, $I$ its ideal such that $A/I$ is SAW*. Suppose that $b\in A/I$, $p(b)=0$, $\|b\|=1$, and a family of projections $\{p_i\}$ associated with $b$ by (\ref{projections})
can be lifted
to a family of
 projections in $A$ with sum 1. Then there is a lift $\tilde b \in A$ of $b$ such that $p(\tilde b)=0, $ $\|\tilde b\|=1$.
\end{lemma}
\begin{proof}
Let $s, c$ be as above and
projections $\{Q_i\}$ in $A$ be lifts of $\{p_i\}$ with sum 1. By (\ref{1}), $p_i(c-\lambda_i)p_i$ is a nilpotent of
order $k_i$ in the $C^*$-algebra $p_i (A/I) p_i$.  By Lemma \ref{aAa}, $p_i (A/I) p_i = Q_i A Q_i/Q_i I Q_i$. It follows from Lemma
\ref{SAW} and Theorem \ref{Olsen, Pedersen} that there are $e_j^{(i)}, d_j^{(i)}\in p_i (A/I) p_i$ and their lifts
$E_j^{(i)}, D_j^{(i)}\in Q_i A Q_i$ satisfying (\ref{2}) and (\ref{3}) such that
$$p_i(c-\lambda_i)p_i =
\sum_{j=1}^{k_i-1}(e_{j-1}^{(i)}-e_j^{(i)})(c-\lambda_i)d_j^{(i)}$$
and, for any lift $C$ of $c$,
$$\sum_{j=1}^{k_i-1}(E_{j-1}^{(i)}-E_j^{(i)})(C-\lambda_i)D_j^{(i)}$$
 is a nilpotent lift of $p_i(c-\lambda_i)p_i$ of order $k_i$.
Let $q:A\to A/I$ be the canonical surjection.
 We have  \begin{multline*} q \left(\sum_{i=1}^n\sum_{j=1}^{k_i-1}\left((
E_{j-1}^{(i)}-E_j^{(i)})( C-\lambda_i) D_j^{(i)} +
\lambda_i Q_i\right)\right) =\\
\sum_{i=1}^n\sum_{j=1}^{k_i-1}\left((
e_{j-1}^{(i)}-e_j^{(i)})(c-\lambda_i)d_j^{(i)} + \lambda_i
p_i\right) =
 \sum_{i=1}^n p_i(c-\lambda_i)p_i + \lambda_i p_i = \sum_{i=1}^n p_icp_i = c.  \end{multline*}
Denote $\sum_{j=1}^{k_i-1}( E_{j-1}^{(i)}-E_j^{(i)})(
C-\lambda_i) D_j^{(i)}$ by $X_i$. Since
$X_i^{k_i}=0$, $X_iQ_j=Q_jX_i=0$ if $i\neq j$, and $\sum Q_i=1$ we
have
\begin{multline*}p\left(\sum_{i=1}^n\sum_{j=1}^{k_i-1}\left((
E_{j-1}^{(i)}-E_j^{(i)})( C-\lambda_i) D_j^{(i)} +
\lambda_i Q_i\right)\right) = \\ \prod_{l=1}^n \left(\sum_{i=1}^n
\left(X_i + \lambda_i Q_i\right) - \lambda_l\right)^{k_l} =
\prod_{l=1}^n
\left(\sum_{i=1}^n\left(X_i+(\lambda_i-\lambda_l)Q_i\right)\right)^{k_l}
= \\ \prod_{l=1}^n \left(X_l+\sum_{i\neq
l}(X_i+\lambda_i-\lambda_l)Q_i\right)^{k_l} = \prod_{l=1}^n
\left(\sum_{i\neq l}(X_i+\lambda_i-\lambda_l)Q_i\right)^{k_l}
=0.\end{multline*}

It follows from (\ref{3}) that $\sum_{j=1}^{k_i-1} (
E_{j-1}^{(i)}-E_j^{(i)}) D_j^{(i)} = Q_i$ whence
$$\sum_{i=1}^n\sum_{j=1}^{k_i-1}\left((
E_{j-1}^{(i)}-E_j^{(i)})( C-\lambda_i) D_j^{(i)} +
\lambda_i Q_i\right) = \sum_{i=1}^n\sum_{j=1}^{k_i-1}(
E_{j-1}^{(i)}-E_j^{(i)}) C D_j^{(i)}. $$

Let $S\in A$ be any invertible positive lift of $s$. Since
$b= s^{-1/2}cs^{1/2}$ the set $$\mathcal E =
\left\{\sum_{i=1}^n\sum_{j=1}^{k_i-1}S^{-1/2}(
E_{j-1}^{(i)}-E_j^{(i)})S^{1/2}XS^{-1/2} D_j^{(i)}S^{1/2}\;|
X\in A,\; q(X) = b\right\}$$ consists of lifts of $b$
and we have $ p(X)=0$, for any $X\in \mathcal E$. Then $\overline{\mathcal E}$ also consists of lifts of $b$ and
$ p(X)=0$, for any $X\in \overline{\mathcal E}$.

Define an elementary operator $T:A \to A$ by $$T(X)
= \sum_{i=1}^n\sum_{j=1}^{k_i-1}S^{-1/2}(
E_{j-1}^{(i)}-E_j^{(i)})S^{1/2}XS^{-1/2} D_j^{(i)}S^{1/2}.$$
Then $\mathcal E = \{TX |X\in A,\; q(X) = b\}$.
By Corollary \ref{new}, in $\overline{\mathcal E}$ there is an element of norm 1.
\end{proof}

\begin{theorem} Let $p(t) ~=
(t-\lambda_1)^{k_1}\ldots (t-\lambda_n)^{k_n}$ be a polynomial such that $k_i\ge 2$, $i=1, \ldots, n$. The universal $C^*$-algebra of a relation $p(b)=0,
\|b\|\le 1$ is semiprojective.
\end{theorem}
\begin{proof} By
Theorem 14.1.7 of \cite{Loring}, it suffices to prove that if $b\in M(I)/I$, where $I = \overline {\bigcup
I_n}$, $I_1\subseteq I_2\subseteq \ldots$, $p(b)=0$, $\|b\|\le 1$, then there exists $N$ and $ \tilde b\in M(I)/I_N$ such that $q_{n, \infty}(
\tilde b) = b$, $p( \tilde b )=0$, $\| \tilde b\|\le 1$. Here by $q_{n, \infty}$ we denote the canonical surjection from $M(I)/I_N$ onto $M(I)/I$.   Clearly we can assume $\|b\|=1$ and look
for $ \tilde b $ with $\| \tilde b\| =1$.

 By \cite{Blackadar}, there exists $N$ such that $\{p_i\}$ can be lifted to a
family $\{Q_i\}$ of projections in $M(I)/I_N$ with sum 1. Denote $M(I)/I$ by $A$, $I/I_N$ --- by $J$. Then $A/J\cong M(I)/I$. By
Kasparov's Technical Theorem, any corona-algebra is SAW*, so the assertion follows now from Lemma \ref{main}.
\end{proof}

Recall that a $C^*$-algebra is residually finite-dimensional (RFD) if the intersection of kernels of its finite-dimensional representations is zero.

\begin{theorem} Let $p(t) ~=
(t-\lambda_1)^{k_1}\ldots (t-\lambda_n)^{k_n}$ be a polynomial such that $k_i\ge 2$, $i=1, \ldots, n$. The universal $C^*$-algebra of a relation $p(b)=0,
\|b\|\le 1$ is residually finite-dimensional.
\end{theorem}
\begin{proof} Let $H$ be $l^2(\mathbb N)$. We will identify $M_n$ with $B(l^2\{1, \ldots, n\})\subseteq B(H)$.
Let $\mathcal A\subseteq \prod M_n$ be a $C^*$-algebra of all sequences $\{x_n\}$ such that $\ast$-strong limit
$x_n$ exists and let $\mathcal I$ be an ideal of all sequences $\{x_n\}$ such that $\ast$-strong limit $x_n$ is
zero. Then a map $i: \mathcal A/\mathcal I \to B(H)$ sending each sequence $\{x_n\} $ to its $\ast$-strong limit
is an isomorphism. The algebra $B(H)$ is SAW*. Indeed if $T, S\in B(H)$, $TS=0$ then for $Z$ equal to the
projection on $(Ran S)^{\perp}$ one has $TZ=T$, $ZS=0$. Thus $\mathcal A/\mathcal I$ is SAW*. Let $p_1, \ldots,
p_n$ be a family of  projections in $\mathcal A/\mathcal I$ with sum 1. Let $H_j$ be a range of projection
$i(p_j)$, $\{e_i^{(j)}\}_{i\in \mathbb N}$ be an orthonormal basis in $H_j$. Then $e_1^{(1)}, \ldots, e_1^{(n)},
e_2^{(1)}, \ldots, e_2^{(n)}, \ldots$ is an orthonormal basis in $H$. Let $R_k$ be a projection on the first k
basis vectors. Then  sequences $\{R_kP_jR_k\}_{k\in \mathbb N} \in \mathcal A$ are projections with sum 1 and
are lifts of $P_j$. By Lemma \ref{main} each $\ast$-homomorphism from the universal $C^*$-algebra of a
relation $p(b)=0, \|b\|\le 1$ to $\mathcal A/\mathcal I$ lifts to an $\ast$-homomorphism to $\mathcal A$. Let
$\pi$ be the universal representation of the universal $C^*$-algebra of a relation $p(b)=0, \|b\|\le 1$.
Then any lift of $i^{-1}\circ \pi$ has to be injective and hence gives a separating family of finite-dimensional
representations.
\end{proof}



\begin{theorem}\label{Calkin} Let $T\in B(H)$ and $p(t) ~=
(t-\lambda_1)^{k_1}\ldots (t-\lambda_n)^{k_n}$ be a polynomial such that $k_i\ge 2$, $i=1, \ldots, n$.  Suppose that
$p(T)$ is compact. Then there is a compact operator $K$ such that $p(T+K)=0$, $\|T+K\|=\|T\|_e$.
\end{theorem}
\begin{proof} Let $b$ be the image of $T$ in Calkin algebra $B(H)/K(H)$, $C= \|T\|_e$. Then $p(b)=0$, $\|b\| = C$.
The assertion follows now from Lemma \ref{main}, Kasparov's Technical Theorem and  the well-known fact that
orthogonal projections (with sum 1) in Calkin algebra can be lifted to orthogonal projections (with sum 1) in
$B(H)$.
\end{proof}

Now we consider the case of a polynomial whose roots are real and have multiplicity one.

\begin{theorem}\label{multiplicity1} Let $p(t) = (t-t_1)\ldots (t-t_n) $ be a polynomial such that $t_i\in \mathbb R$, $t_i\neq t_j$ when $i\neq j$. Then the universal $C^*$-algebra of the relations $p(x)=0$, $\|x\|\le C$ is semiprojective.
\end{theorem}
\begin{proof}   Let $A$ be a $C^*$-algebra, $\{I_n\}$  an increasing chain of ideals of $A$ and  $I=\overline{\bigcup I_n}$. By  $q_N:A\to A/I_N$ we will denote the canonical surjection and for an element $a$ in any $A/I_N$, $N\in \mathbb N$, by $\dot a$ we denote its image in $A/I$.

 Suppose $x\in A/I$, $p(x)=0$, $\|x\|\le C$. We have to show that $x$ has a lift $X$ in $A/I_N$, for sufficiently large $N$, such that $p(X)=0$, $\|X\|\le C$.

Let $s$ and $p_i$ be as in (\ref{projections}) and let $y= \sum_{i=1}^n t_i p_i$. By (\ref{1}) and (\ref{4}), $x = s^{-1/2}ys^{1/2}$.
 Fix any invertible lift $S$ of $s$.   Define $T: A\to A$ by $$T(a) = S^{-1/2}aS^{1/2},$$ $a\in A$. Let $Y_0$ be any selfadjoint lift of $y$.  We have $\|\dot{(TY_0)}\| = \|s^{-1/2}ys^{1/2}\| = \|x\| \le C$. Let $$E = \{T(Y_0) + T(i)\;| \; i\in I_{sa}\}.$$
 Let $\{e_{\lambda}\}$ be a quasicentral approximate unit for $I\triangleleft A$ and $$i_{\lambda}= 1-(1-e_{\lambda})^{1/2}.$$ Then $\{i_{\lambda}\}$ is also a quasicentral approximate unit and we have  \begin{multline*} dist (T(Y_0), TI_{sa}) =  \inf_{e\in E} \|e\| \le \limsup \|T(Y_0) - T\left(i_{\lambda}Y_0+Y_0i_{\lambda} - i_{\lambda}Y_0i_{\lambda}\right)\| = \\ \limsup \|T\left((1-i_{\lambda})Y_0(1-i_{\lambda})\right)\|
     = \limsup \|(T(Y_0))(1-i_{\lambda})^2\| = \\ \limsup \|(T(Y_0))(1-e_{\lambda})\| = \|\dot{T(Y_0)}\| \le C.\end{multline*}   Since $S$ is invertible, $T$ is also invertible and hence
 $TI_{sa}$ is a closed subspace of $TA_{sa}$.   By Theorem \ref{prox}, there is $i_0\in I_{sa}$ such that $dist (TY_0, TI_{sa}) = \|TY_0-Ti_0\|$. Let $$Y = Y_0-i_0.$$ Thus we found a self-adjoint lift $Y$ of  $y$ such that \begin{equation*} \|S^{-1/2}YS^{1/2}\|\le C.\end{equation*}

Let $O_i$, $i=1, \ldots, n$, be neighborhoods of $t_1, \ldots, t_n$  such that $O_i\bigcap O_j=\varnothing$, when $i\neq j$. There is
$\delta>0$ such that if $|p(t)|\le \delta$ then $t\in \bigcup O_i$.
Since $p(y)=0$, by [Prop. 2.14, \cite{Blackadar}] for sufficiently large N, $\|p(q_N(Y))\|\le \delta$   and hence
     $$\sigma(q_N(Y)) \subseteq \bigcup_i O_i.$$  Let $D_C$ be the disk of radius $C$. Let $f$ be a continuous function on $D_C$ such that
\begin{itemize}
\item $f(t_i)=t_i$ when $t\in O_i$;
\item $|f(t)|\le C$, for any $t\in D_C$.
\end{itemize}
Let $X=q_N(S^{-1/2}f(Y)S^{1/2}) = q_N(S^{-1/2})f(q_N(Y))q_N(S^{1/2})$. Then $\dot X = x$. Since  $f(q_N(Y))$ is normal and  $\sigma(f(q_N(Y))) = \{t_1, \ldots, t_n\}$, we have $$(X-t_1)\ldots (X-t_n) = q_N(S^{-1/2})(f(q_N(Y))-t_1)\ldots (f(q_N(Y))-t_n)q_N(S^{1/2})=0.$$

Let $p_k$ be polynomials that converge to $f$ uniformly on $D_C$  such that $$\sup_{t\in D_C}|p_n(t)|\le C.$$
 Then \begin{multline*} \|X\| = \|q_N(S^{-1/2}f(Y)S^{1/2})\| = \lim_{k\to \infty} \|q_N(S^{-1/2})p_k(q_N(Y))q_N(S^{1/2})\| = \\ \lim_{k\to \infty} \|p_k(q_N(S^{-1/2}YS^{1/2}))\| \end{multline*} and by the von Neumann inequality $$\|X\| \le \lim_{k\to \infty} \sup_{|t|\le \|q_N(S^{-1/2}YS^{1/2})\|} |p_k(t)| \le \lim_{k\to \infty} \sup_{t\in D_C} |p_k(t)|\le C.$$

\end{proof}

\begin{theorem}\label{Calkin2} Let $p(t) = (t-t_1)\ldots (t-t_n) $ be a polynomial such that $t_i\in \mathbb R$, $t_i\neq t_j$ when $i\neq j$ and let  $T\in B(H)$.  Suppose that
$p(T)$ is compact. Then there is a compact operator $K$ such that $p(T+K)=0$, $\|T+K\|=\|T\|_e$.
\end{theorem}
\begin{proof} Let $x$ be the image of $T$ in Calkin algebra. Let $C= \|T\|_e$. It is sufficient to find a pre-image $X\in B(H)$ of $x$ such that $p(X)=0$ and $\|X\|=C$. As in the proof of Theorem \ref{multiplicity1} we define $s$, $y$ and $S$ and find a selfadjoint pre-image $Y$ of $y$ such that
\begin{equation*} \|S^{-1/2}YS^{1/2}\|\le C.\end{equation*}
Since $p(Y)$ is compact, its spectrum is either finite or a sequence vanishing at $0$. Hence the spectrum of $Y$ is either finite or countable with accumulation points $t_i$, $i=1, \ldots, n$. In either case one can find neighborhoods $O_1, \ldots, O_n$ of $t_1, \ldots, t_n$  respectively,  such that $O_i\bigcap O_j=\varnothing$, when $i\neq j$, and $\sigma(Y)\subset \bigcup O_i$. Now it remains to repeat the arguments from the proof of Theorem \ref{multiplicity1}.
\end{proof}

 \end{document}